\documentclass[12pt]{amsart} 
\usepackage[all]{xy}

\theoremstyle{remark}
\newtheorem{rem}{Remark}

\theoremstyle{definition}
\newtheorem{defi}{Definition}

\newcommand{\he}{\simeq}                

\theoremstyle{theorem}
\newtheorem{theo}{Theorem}[section]
\newtheorem{cor}[theo]{Corollary}
\newtheorem{prop}[theo]{Proposition}
\newtheorem{lem}[theo]{Lemma}

\newcommand{\bib}[6]{ \bibitem{#1}  #2, \textit{ #3\/}, #4 \textbf{#5} #6}

\begin{document}

\title{Homotopical Dynamics: Suspension and Duality}
\author{Octavian Cornea
}\address{ 
Universit\'{e} de Lille 1\\
U.F.R. de Math\'{e}matiques \& U.R.A. 751 au CNRS\\
59655 Villeneuve D'Ascq, France}

\email{cornea@gat.univ-lille1.fr}
\urladdr{http://www-gat.univ-lille1.fr/\~{}cornea/octav.html}
\date{December 1997}

\begin{abstract}
Flow type suspension and homotopy suspension agree for attractor-repellor 
homotopy data.
 The connection maps associated in Conley index theory to an attractor-repellor 
decomposition
with respect to the direct flow and its inverse are Spanier-Whitehead duals in the
stably parallelizable context and are duals modulo a certain Thom construction in
general.
\end{abstract}

\maketitle

\section{Introduction.}
The Morse inequalities as well as their
degenerate counterpart, the Lusternik-Schnirelmann inequality \cite{LusternikS}, are early
 examples of applications of homotopy theory to the study of 
smooth flows. 
More recently, the subject has developed along three directions:

 -Various numerical homotopical 
invariants related to the Lusternik-\\Schnirelmann category have 
been investigated  intensively.   

-Particular attention has been given to Morse-Smale flows 
 where the use of homotopical methods has had much success.  

 -The index theory of Conley \cite{Conley} has introduced homotopical methods in
the analysis of general flows.

The last few years have seen a certain convergence of these topics having at the 
center the notion of stabilization.
For example, given a $2$-connected manifold $M$, for $k$
large enough, one can construct functions on the manifold $M\times D^{k}$  
pointing inwards on the boundary and with no more than $cat(M)+2$ critical 
points \cite{Cornea1}. Attaining the best possible lower bound, $cat(M)+1$, depends, for now, 
on showing that the L.S.-category and another homotopical invariant,
 the cone-length of $M$ agree \cite{Cornea4}. This seems likely for most
closed manifolds \cite{CorneaFL} but for arbitrary $CW$-complexes it is not true and recently
 have been constructed \cite{Dupont}
new homotopical examples for which the two invariants differ. 

Conversely, one can look for critical point estimates for
functions that are defined on $M\times D^{k}$ with $M$ closed and which are
pointing out on $M\times D^{k-t-1}\times S^{t}$ and in on the complement
of this set in $\partial (M\times D^{k})$.  The lower bound in this case is
$cat(M,M\times S^{t})$ \cite{Cornea3}. 
For any $CW$-complex, this invariant can be shown \cite{Moyaux} to be strictly bigger
 than a more stable version of the L. S.-category $\sigma^{(t+1)} - cat$ 
 (defined as the least $n$ for which the $(t+1)$- suspension of the $n$-th filtration in Milnor's 
classifying construction applied to $\Omega M$ admits a homotopy section \cite{Vandembroucq}). 

This class of functions contains the functions quadratic at infinity
studied in relation to the Arnold conjecture in symplectic geometry. 
They are at the center of the "stable" Morse and L.S. theory of
Eliashberg and Gromov \cite{EliashbergG}. Because of the Arnold conjecture
one does then expect that $cat(M,M\times S^{t})\geq cat(M)+1$. 
This inequality is known to hold in many cases and it
 implies \cite{Cornea3} $cat(M\times S^{t})=
cat(M)+1$. In its turn, the validity of this equality for all $CW$-complexes is known as 
the Ganea conjecture. 
The conjecture has been recently
disproved \cite{Iwase}. However, for symplectic, closed manifolds of symplectic
form $\omega$ and with $\omega|_{\pi_{2}(M)}=0$
 by \cite{OpreaR} one has  $\sigma^{(t+1)} -cat(M)=cat(M)=
cat(M,M\times S^{t})-1=dim(M)$.

Finally,  again in relation to the Arnold conjecture and Floer homology,
 Cohen, Jones and Segal \cite{CohenJS} have introduced a stable homotopy point of 
view in studying Morse-Smale flows (see also \cite{Hurtubise}).\medskip

All of this suggests that in the flow context suspending should be understood to 
mean taking the product of the given flow with 
the gradient flow of a quadratic form on $\mathbf{R^{n}}$. 

The present paper is concerned with the next step in the study of
the homotopical properties of flows: attractor-repellor data. 
We address two main issues in this context.

\begin{itemize}
\item Understand the behaviour of attractor-repellor data under suspension. 
\item Determine in what measure and in what sense the homotopical information provided
by a smooth flow and its inverse are redundant.
\end{itemize}

The key point is that suspension is central in dealing with the second matter.

Attractor-repellor data is best understood for Morse-Smale functions.
Passage through a critical point of a Morse-Smale function corresponds to
the attachement of a cell. Hence, for two consecutive critical points one has a relative
attaching map $d:S^{q-1}\longrightarrow S^{k}$ where $k$ is the index of the first
point and $q$ that of the second. This map has been described by Franks in \cite{Franks}.
The negative of the original function being also Morse-Smale there is also
a relative attaching map $d':S^{n-k-1}\longrightarrow S^{n-q}$ with
$n$ being the dimension of the supporting manifold $M$.  Franks has also shown that,
if $M$ is stably parallelizable, then $d$ and $d'$ are stably the same up to sign.

Conley index theory allows one to define analogues of these maps
for general attractor-repellor pairs of some general flow $\gamma$ (which will
be supposed here to be defined on a smooth manifold and to be itself smooth). 
Assume $S$ is an isolated invariant set of $\gamma$ and $(A, A^{\ast})$
is an attractor-repellor decomposition of $S$. Then one has the map
$\delta: c_{\gamma}(A^{\ast})\longrightarrow \Sigma c_{\gamma}(A)$
called connection map of the attractor-repellor pair $(A, A^{\ast})$. 
Let $-\gamma$ be the inverse flow of $\gamma$,
$-\gamma_{t}(x)=\gamma_{-t}(x)$. As $(A^{\ast},A)$ is an attractor-repellor
decomposition of $S$ with respect to $-\gamma$ we also have the inverse connection map
 $\delta':c_{-\gamma}(A)\longrightarrow \Sigma c_{-\gamma}(A^{\ast})$.
 
Here $c_{\gamma}(-)$ is the Conley index of the respective
invariant set. Specialized to the Morse-Smale case
and consecutive critical points $\delta=\Sigma d$ and $\delta'=\Sigma d'$.\medskip

The results of the paper imply two main things:
\begin{itemize}
\item Flow type suspension and homotopy suspension agree for attractor-repellor 
homotopy data.
\item The maps $\delta$ and $\delta'$ are Spanier-Whitehead duals in the
stably parallelizable context and are duals modulo a certain Thom construction in
general.
\end{itemize}
The duality result implies
the (co)-homological results of McCord \cite{McCord}. In the
stably parallelizable case the duality together with the suspension result imply
 that, up to flow suspension, the direct and inverse attractor-repellor homotopy data 
are redundant. 
Conley index theory works quite well for continuous flows on 
metrizable spaces, hence in singular or stratified contexts and parts
of the theory extend even to semiflows. Our
suspension result is likely to extend also to the most general setting. Not so the
duality. This result seems to be specific to flows on manifolds and shows that 
Spanier-Whitehead duality is a key ingredient in understanding the related homotopical
data. Possibly,
the smootheness restriction might be relaxed.

Of course, the result of Franks on the stable equality up to sign of $d$ and $d'$ 
is also extended (as Spanier-Whitehead duality of maps between spheres is stable equality
up to sign). His paper \cite{Franks} is closest in spirit to ours
 even if, due to the fact that for general flows one has no control 
on stable and unstable manifolds, our methods are forced to be quite different.

In trying to make the paper more easely accessible we have included, in the second section,
 a relatively detailed recall of Spanier-Whitehead duality and of the main properties
of the Conley index.
 
The third section contains the precise statements and proofs and
is organized as follows. We first discuss the Thom construction in the setting of the
Conley index. This is needed to prove the suspension result and it also
reduces proving the duality to the study of a flow defined on a sphere. For such a flow
the duality of the two connection maps is obtained by identifying each of them
to a connectant in a certain cofibration sequence and then showing that
 these two cofibration sequences are dual.

The paper is written in the language of connected simple systems for two reasons.
Firstly, the constructions related to the Conley index contain many choices and in this way
one is insured of the invariance of the result. Secondly, connected simple systems reflect
precisely the homotopical information that can be extracted from the flow.

\section{Recalls.}

We recall a number of well-known facts, first on Spanier-Whitehead duality
and then concerning the Conley index.

\subsection{Spanier-Whitehead duality.}

\subsubsection{Definitions.}
Two pointed, connected
$CW$-complexes
$X$, $X'$ are Spanier-Whitehead dual  \cite{SpanierW},\cite{Spanier} 
if there is a map $\mu:X\wedge
X'\longrightarrow S^{m}$ such that the slant product
$\mu^{\ast}(i_{m})^{\ast}/ : H_{q}(X')\longrightarrow H^{m-q}(X)$ gives an
isomorphism.
 The map $\mu$ is called an $m$-duality. 
For completeness recall (see for example \cite{Greenberg}) that the slant 
product $\alpha/\gamma$ with
$\alpha\in H^{p+q}(A\times B)$ and $\gamma\in H_{p}(A)$ is a
cohomology class in $H^{q}(B)$ that satisfies $<\alpha/\gamma, \beta>=
<\alpha,\gamma\times\beta>$ for all $\beta\in H_{q}(B)$,
  where $\times$ is the exterior cross product
and $< , >$ is the usual cohomology-homology pairing. The natural map
$A\times B\longrightarrow A\wedge B$ allows one to use in the slant product
a class $\alpha\in H^{p+q}(A\wedge B)$ by pulling it back to the product $A\times B$. 

Spanier-Whitehead duality behaves well under
suspension: if $X$ and $X'$ are
$m$-duals with duality map $\mu$, then
$\Sigma^{k}X$ and $\Sigma^{q}X'$ are $m+k+q$-duals with respect to
the obvious suspension of $\mu$. 
 Assume $X$, $X'$
are duals by an $m$-duality, $\mu$, and $Y$, $Y'$ are duals by
an $m$-duality $\nu$. Two pointed maps $f:X\longrightarrow Y$,
$g:Y'\longrightarrow X'$ are duals if for some 
$k$,$k'$ the following diagram homotopy commutes:
$$\xy\xymatrix@-2pt{
 \Sigma^{k}X\wedge\Sigma^{k'} Y' \ar[rr]^{f\wedge
1}\ar[d]_{1\wedge g}& &\Sigma^{k}Y\wedge
\Sigma^{k'}Y' \ar[d]^{\nu}\\
\Sigma^{k}X\wedge
\Sigma^{k'}X' \ar[rr]_{\mu}& &S^{m+k+k'} }
\endxy $$
In general, only by allowing $k$ and
$k'$ to be sufficiently large this notion of duality gives a duality isomorphism between 
(stable) classes of maps $D_{m}(\mu,\nu):\{X,Y\}\approx\{Y',X'\}$.

Spanier-Whitehead duality implies Poincar\'e duality via the Thom isomorphism. 
For example, let $(W^{n};V_{0},V_{1})$ be  a smooth manifold triad with 
$V_{0}\coprod V_{1}=\partial W$. Let $\nu$ be the stable normal bundle of $W$ and assume
it is of rank $m>>n$. Then  
$W/V_{0}$ is $(m+n)$- Spanier-Whitehead dual to $T^{\nu}(W)/T^{\nu}(V_{1})$ where
$T^{\nu}(-)$ is the respective Thom space \cite{Atiyah}
(recall that the Thom space of an orthogonal bundle $\eta$ over a space 
$X$ is obtained form the total space of the unit disk bundle associated to $\eta$ by identifying
to a point the total space of the unit sphere bundle). We have the Thom
isomorphism $H^{m+k}(T^{\nu}(W)/T^{\nu}(V_{1}))\simeq H^{k}(W/V_{1})$.
On the other hand we also have the Spanier-Whitehead duality isomorphism
$H^{m+k}(T^{\nu}(W)/T^{\nu}(V_{1}))\simeq H_{n-k}(W/V_{0})$.

\subsubsection{Construction of duality maps.}
Here is one way to obtain duality maps \cite{SpanierW}. Assume 
$A\subset S^{m+1}$, $B\subset S^{m+1}$ such that the inclusion
induces an isomorphism $H_{\ast}(B)\longrightarrow H_{\ast}(S^{m+1}-A)$,
and $A$, $B$ are connected and disjoint.  Define a map in the following way :  
fix a point $p\in S^{m+1}-(A\bigcup B)$ and identify
$S^{m+1}-\{p\}$ with ${\bf R}^{m+1}$. We then define $v:A\times
B\longrightarrow S^{m}$ by $v(x,y)=(x-y)/||x-y||$. Connectivity implies
that this map is null-homotopic when restricted to $A\vee B$. It defines
a map $A\wedge B\longrightarrow S^{m}$ which is a duality map.
One can identify dual maps in this setting:
if $A'$, $B'$ is another pair of subspaces of $S^{m+1}$ that satisfies the
conditions above and if $A\stackrel{i}{\longrightarrow}A'$ and 
$B'\stackrel{i'}{\longrightarrow}B$ are inclusions, then $i$ and $i'$
are dual with respect to the duality maps constructed as above.

\subsection{ The Conley index.}

We will use here and in the rest of the paper the conventions of \cite{Salamon}.

\subsubsection{Connected simple systems.}
 Recall first that 
the pointed homotopy category of pointed topological spaces $Ho-{\mathbb T}$
 has as objects pointed topological spaces and as morphisms homotopy classes of 
pointed (continuous) maps. A connected simple system is a subcategory ${\mathbb C}$ of
$Ho-{\mathbb T}$ such that for any two objects $X_{1}$, $X_{2}$ of ${\mathbb C}$  the set 
$mor_{\mathbb C}(X_{1},X_{2})$ consists of exactly one element. Any map 
$X_{1}\longrightarrow X_{2}$ whose homotopy class belongs to
$mor_{\mathbb C}(X_{1},X_{2})$  is called a comparison
map. Clearly, all comparison maps are homotopy equivalences.
A morphism between two connected simple systems $X$ and $Y$ 
associates to each pair $X_{0}\in Ob(X)$ and
$Y_{0}\in Ob(Y)$ a unique homotopy class of maps $X_{0}\longrightarrow Y_{0}$
with the obivous compatibility properties. We will sometimes call such a morphism
a \textit{map} of connected simple systems.
Connected simple systems form
also a category ${\mathbb CS}$ (see also \cite{McCord2}).  Let $X, Y\in Ob({\mathbb CS})$. 

Assume
that we have a (non-void)
 set ${\mathbb P}$ of pairs $(A,B)$ with $A\in Ob(X)$ and
$B\in Ob(Y)$ and maps $f_{A,B}:A\longrightarrow B$ that are defined for each such pair. 
We say that these maps
\textit{induce} a morphism $f:X\longrightarrow Y$ if
for any two pairs $(A,B)$ and $(A',B')$ in ${\mathbb P}$ the maps $f_{A,B}$ and
$f_{A',B'}$ commute, up to homotopy, with the obvious comparison maps and $f$
is the unique morphism of connected simple systems represented by the maps $f_{A,B}$.
 Of course, defining $f$ is possible even if ${\mathbb P}$ has a single element.
However, in many of the constructions below the maps
$f_{A,B}$  are naturally defined for many pairs $(A, B)$ and it is useful
to know that the induced morphism does not depend on the particular
choice of the pair. Hence, in the following we will speak of an \textit{induced map}
at the level of connected simple systems by always considering ${\mathbb P}$
to be as large as is allowed by the definition of the maps $f_{A,B}$.  

\subsubsection{Cofibration sequences.}
An obvious functor $\Sigma : {\mathbb CS}\longrightarrow
{\mathbb CS}$ is defined by taking the supension of spaces and maps. By definition,
a cofibration sequence of connected simple systems is a
triple of maps $X\stackrel{i}{\longrightarrow}Y
\stackrel{p}{\longrightarrow}Z\stackrel{\delta}{\longrightarrow}\Sigma X$
with $X, Y, Z\in Ob({\mathbb CS})$ with the property that there are $X_{0}\in Ob(X)$,
$Y_{0}\in Ob(Y)$, $Z_{0}\in Ob(Z)$, a cofibration sequence
 $X_{0}\stackrel{i_{0}}{\longrightarrow}Y_{0}
\stackrel{p_{0}}{\longrightarrow}Z_{0}$ 
with $[i_{0}]=i$, $[p_{0}]=p$ and if $\delta_{1} :Z_{0}\longrightarrow \Sigma X_{0}$
is the obvious connectant, then $\delta =[\delta_{1}]$. The map $\delta$  is
called the connection map of the cofibration sequence. The cofibration
sequence $X_{0}\longrightarrow Y_{0}\longrightarrow Z_{0}$
is called a representing cofibration sequence.
The cofibration sequence is determined by the choice of $i_{0}$ 
(however,  because push outs in $Ho-{\mathbb T}$ 
do not coincide, in general, with homotopy push outs, the knowledge of $i$
only does not suffice).
For completeness we recall the definition of a cofibration sequence. An inclusion
$j:A\longrightarrow B$ is a cofibration if it satifies the homotopy extension property.
This means that given a homotopy $h:A\times [0,1]\longrightarrow X$ and a map
$F:B\longrightarrow X$ such that $F\circ j=h_{0}$ then there is a homotopy
$H:B\times [0,1]\longrightarrow X$ that extends $h$ and with $H_{0}=F$. This concept
is useful here because if $j$ is a cofibration then the homotopy type of the quotient
$B/A$ is the same as that of $B\bigcup_{j} CA$ ($CA$ is the cone over $A$). 
In view of this, we will sometimes identify $C$ and $B/A$. The
pair of maps $A\stackrel{j}{\longrightarrow} B\stackrel{r}{\longrightarrow} C$ 
form a cofibration sequence if $j$ is a cofibration,
$C=B/A$ and $r$ is the obvious collapsing map.
A cofibration sequence extends  to the right giving
a sequence of maps: 
$A\stackrel{j}{\longrightarrow}B\stackrel{r}{\longrightarrow}
C\stackrel{\delta_{1}}{\longrightarrow}
\Sigma A\stackrel{\Sigma j}{\longrightarrow}\Sigma B\stackrel{\Sigma r}
{\longrightarrow}\cdots$
where $\delta_{1}$ is the collapsing $C \simeq B\bigcup CA\longrightarrow (B\bigcup CA)/B=
\Sigma A$.

\subsubsection{The Conley index.}
We recall now a few elements of Conley index theory \cite{Conley},\cite{Salamon}
for the case of a continuous flow $\gamma: N\times {\bf R}\longrightarrow N$
, $N$ being a compact metric space.
Let $S\subset N$ be an isolated invariant set. A pair $(N_{1},N_{0})$ of compact sets in
$N$ is an index pair for $S$ in $N$ if $N_{0}\subset N_{1}$, $N_{1}-N_{0}$ is a 
neighborhood of
$S$, $S$ is the maximal invariant set in the closure of $N_{1}-N_{0}$, 
$N_{0}$ is positively invariant
in $N_{1}$  and, if for $x\in N_{1}$ there is some $t\geq 0$ such that $\gamma_{t}(x)\not\in N_{1}$,
then there exists a $\tau >0$ with $\gamma_{t}(x)\in N_{1}$ for $0\leq t\leq \tau$ and
$\gamma_{\tau}(x)\in N_{0}$, \cite{Salamon}.

The basic result in this theory is the existence of index pairs inside any isolating
neighborhood. 
 Moreover, for any two index pairs of $S$,  $(N_{1},N_{0})$ and 
$(N'_{1},N_{0})$ there are flow induced comparison map $N_{1}/N_{0}\longrightarrow
N'_{1}/N'_{0}$. In the following we need a uniform choice for the comparison maps.
We will refer to the maps defined by the formula of Lemma 4.7 of \cite{Salamon} 
as the \textit{standard comparison
maps}. This gives the set of quotients $N_{1}/N_{0}$ the structure of a connected
simple system denoted by $c_{\gamma}(S)$ and called the Conley index of $S$ with
respect to $\gamma$.  

Assume that $S$ admits an attractor-repellor decomposition, denoted $(S;A, A^{\ast})$,
 of attractor $A$ and repellor $A^{\ast}$ (this means that if $x\in S-(A\bigcup A^{\ast})$
 then $\lim_{t\rightarrow\infty}\gamma_{t}(x)\in A$ and 
$\lim_{t\rightarrow -\infty}-\gamma_{t}(x)\in A^{\ast}$). 
We have a cofibration sequence of connected simple systems denoted
 by $c_{\gamma}(S;A,A^{\ast})$: 
$c_{\gamma}(A)\stackrel{i}{\longrightarrow}
 c_{\gamma}(S)\stackrel{p}{\longrightarrow} c_{\gamma}(A^{\ast})
\stackrel{\delta}{\longrightarrow}\Sigma c_{\gamma}(A)$
such that $\delta$ as well as $i$ and $p$ are flow defined \cite{Salamon}. In this case $\delta$ is
called the connection map of the given attractor-repellor decomposition.

Moreover, for any triple $N_{0}\subset N_{1}\subset N_{2}$ with the inclusions being
inclusions of NDR's and such that $(N_{2},N_{1})$ is an index pair of $ A^{\ast}$, 
$(N_{1},N_{0})$ is an index pair of $A$ and $(N_{2},N_{0})$ an index pair of $S$
the cofibration sequence 
$N_{1}/N_{0}\longrightarrow N_{2}/N_{0}\longrightarrow N_{2}/N_{1}$
 represents $c_{\gamma}(S;A,A^{\ast})$. 

For $M\subset N$ we use the notation $I(M)$ for the maximal invariant set
inside $M$.

\subsubsection{Duality and connected simple systems.}
The definition of Spanier-Whitehead duality induces one for simple connected systems.
Consider the connected simple system formed by a single space $S^{m}$
and let $X,X'\in Ob({\mathbb CS})$. We have $X\wedge X'\in Ob({\mathbb CS})$.
The two connected simple systems are Spanier-Whitehead duals if there is a morphism
$X\wedge X'\longrightarrow S^{m}$ giving a duality map
whenever restricted to $X_{0}\wedge X'_{0}$, $X_{0}\in Ob(X)$, $X'_{0}\in Ob(X')$ .
A similar definition applies to maps. 

Let $X\stackrel{i}{\longrightarrow}Y\stackrel{p}{\longrightarrow}
 Z\stackrel{\delta}{\longrightarrow}\Sigma X$
and $Z'\stackrel{p'}{\longrightarrow}Y'\stackrel{i'}{\longrightarrow}
 X'\stackrel{\delta'}{\longrightarrow}\Sigma Z'$ be two cofibration sequences
of connected simple systems. They are called dual cofibration sequences
if the pairs of maps $(i,i')$, $(p,p')$ and $(\delta,\delta')$ are respectively dual.
To insure this duality it is enough to find cofibration sequences of spaces:
$X_{0}\stackrel{i_{0}}{\longrightarrow}Y_{0}\stackrel{p_{0}}{\longrightarrow}Z_{0}$
and $Z'_{0}\stackrel{p'_{0}}{\longrightarrow}Y'_{0}\stackrel{i'_{0}}{\longrightarrow}
 X'_{0}$ representing respectively the two above and with $(i_{0},i'_{0})$
and $(p_{0},p'_{0})$ respectively duals. Indeed, in this case, by the basic properties
of Spanier-Whitehead duality \cite{Spanier}, we have that the associated connectant maps
$\delta_{0}:Z_{0}\longrightarrow \Sigma X_{0}$ and
$\delta'_{0}:X'_{0}\longrightarrow \Sigma Z'_{0}$ are Spanier-Whitehead duals.
Of course, under these circumstances, one wants also to show that the duality relation 
obtained does not depend on the choice of
representing cofibration sequences.

\section{Bundles, suspension and duality.}
Let $\psi:F\longrightarrow E\stackrel{p}{\longrightarrow} B$ be a locally trivial
fibration of manifolds (possibly with boundary) with $F$ and $B$ compact, 
and let $\gamma$ be a flow on $B$. A 
lift of $\gamma$ to $E$ is
a flow $\gamma'$ on $E$ such that $p(\gamma'_{t}(x))=\gamma_{t}(p(x))$
for all $x\in E$ and $t\in \mathbf{R}$.

\subsection{The Thom index.} 
We fix a first construction with the next immediate fact.

\begin{lem}Let $S$ be an isolated invariant set of the flow 
$\gamma :B\times \mathbf{R}\longrightarrow B$ and let $\gamma'$ be a lift of $\gamma$
to $E$. Then, $S'=p^{-1}(S)$ is an isolated, invariant set of $\gamma'$ and
there is an induced morphism $c_{\gamma'}(S')\longrightarrow c_{\gamma}(S)$
depending only on $\psi$ and $\gamma$.
This morphism is natural with respect to attractor-repellor pairs.
\end{lem}

\begin{proof} 
Let $(N_{1},N_{0})$ be an index pair of $S$. We first notice that
$S'=I(f^{-1}(\overline{N_{1}-N_{0}}))$. Moreover,
the pair $(f^{-1}(N_{1}),f^{-1}(N_{0}))$ is an index pair for $S'$  (with respect
to the flow $\gamma'$, of course).
Denote $N_{1}'=f^{-1}(N_{1})$ and $N_{0}'=f^{-1}(N_{0})$.
There is a map of pairs $(N_{1}',N_{0}')\longrightarrow (N_{1},N_{0})$
which induces a map of quotients: $p_{N}:N_{1}'/N_{0}'\longrightarrow N_{1}/N_{0}$.
Notice that this map induces one of connected
simple systems in the following way. 

First, let $(K_{1},K_{0})$ be another
index pair of $S$ and let, as above, $(K_{1}',K_{0}')$ be the index pair
of $S'$ with $K_{i}'=p^{-1}(K_{i})$, $i=0,1$. There are standard, 
flow defined \cite{Salamon}, comparison
maps, $N_{1}/N_{0}\longrightarrow K_{1}/K_{0}$ and
$N_{1}'/N_{0}'\longrightarrow K_{1}'/K_{0}'$. By construction, these maps
commute with the projections $p_{N}$ and $p_{K}$.

Now assume that $(N_{1}'',N_{0}'')$ is another
index pair of $S'$ (and not necessarily one obtained as a preimage of one of
$S$) and consider also an index pair of $S$, $(T_{1},T_{0})$.
We define a map $N_{1}''/N_{0}''\longrightarrow T_{1}/T_{0}$ as the 
composite $N_{1}''/N_{0}''\longrightarrow N_{1}'/N_{0}'\stackrel{p_{N}}{\longrightarrow}
N_{1}/N_{0}\longrightarrow T_{1}/T_{0}$ where the first and last maps
in this composition are the standard comparison maps.
The homotopy class of this map does not depend on the choice of the 
pair $(N_{1},N_{0})$. 

To prove naturality with respect to attractor-repellor sequences consider 
an attractor-reppellor decomposition $(S;A,A^{\ast})$ for an invariant set $S$ and
 a triple $N_{0}\subset N_{1}\subset N_{2}$ such that the
pair $(N_{1},N_{0})$  is an index pair for $A$, the pair
$(N_{2}, N_{1})$ is an index pair for $A^{\ast}$ and $(N_{2},N_{0})$ is an index
pair for $S$. Following the notations above we denote by $'$ the respective
 preimages by $p$ of these sets. It is clear that $A'$ and $(A^{\ast})'$ 
form an attractor-repellor pair. The result now follows from the
commutativity of the diagram:
$$ \xy\xymatrix@-2pt{
N_{1}'/N_{0}' \ar[r]\ar[d]& N_{2}'/N_{0}'\ar[d]\ar[r]&N_{2}'/N_{1}'\ar[d]\\
N_{1}/N_{0}\ar[r]&N_{2}/N_{0} \ar[r]&N_{2}/N_{1} }

\endxy
$$
The vertical arrows are given, as above,  by projection. The rows
are the cofibration sequences that represent
the attractor-repellor cofibration sequences of connected-simple systems. Moreover, note that
the right vertical map up is induced
on cofibres by the left and middle projections and this implies that, by extending the 
cofibrations to the right, we will continue to get commutative squares.
\end{proof}

For a pointed space $X$ let, as before, $CX$ (or $C(X)$) be the (reduced) cone on $X$.

\begin{defi} With the notations in the proof of Lemma 3.1, the Thom index of 
$S$ with respect to $\psi$ and $\gamma$ is the connected simple system
$\overline{c}_{\gamma}^{\psi}(S)$ having as objects the spaces $(N_{1}/N_{0})
\bigcup_{p_{N}}C(N'_{1}/N'_{0})$ and with the morphisms given by the
homotopy classes of the maps\\  $(N_{1}/N_{0})
\bigcup_{p_{N}}C(N'_{1}/N'_{0})\longrightarrow (K_{1}/K_{0})
\bigcup_{p_{K}}C(K'_{1}/K'_{0})$ induced by the standard comparison maps. 
\end{defi}

\begin{lem}In the context of the lemma above and assuming that
$(S;A, A^{\ast})$ is an attractor-repellor decomposition there is a
cofibration sequence of connected simple systems:
$\overline{c}_{\gamma}^{\psi}(A)\longrightarrow 
\overline{c}_{\gamma}^{\psi}(S)\longrightarrow 
\overline{c}_{\gamma}^{\psi}(A^{\ast})\longrightarrow 
\Sigma\overline{c}_{\gamma}^{\psi}(A)$.
\end{lem}

\begin{proof} The naturality given by the previous lemma implies the existence of the
maps in the statement. These maps give a cofibration sequence because
of the general fact that if in the diagram:
$$\xy\xymatrix@-2pt{
A'\ar[d]\ar[r]&B'\ar[d]\ar[r]&C'\ar[d] \\
A\ar[d]\ar[r]&B \ar[d]\ar[r]&C\ar[d] \\
A''\ar[r]&B''\ar[r]&C'' }
\endxy $$
the top two horizontal rows as well as the columns are cofibration sequences,
then the third horizontal row is also a cofibration sequence.
In our case, the top row contains the spaces representing the Conley index
of $\gamma'$ . The spaces in the middle row are obtained from the spaces giving
the Conley indexes with respect to $\gamma$ by pasting the reduced cylinders 
of the projections. The vertical maps are inclusions in the free ends of these cylinders.
\end{proof}

Assume that $(S;A, A^{\ast})$ is an attractor-repellor decomposition for the
flow $\gamma$. Recall that we denote by $c_{\gamma}(S;A,A^{\ast})$
the attractor-repellor cofibration sequence associated to this decomposition. We will
denote by $\overline{c}_{\gamma}^{\psi}(S;A,A^{\ast})$ the cofibration
sequence given by Lemma 3.2.

\begin{rem}
Notice that this also applies to the limit case when the total space of $\psi$
is void. In this case, $\overline{c}_{\gamma}^{\psi}(S;A,A^{\ast})$ is isomorphic to
$c_{\gamma}(S;A,A^{\ast})$ 
(recall that $X/\emptyset$ is the disjoint union of $X$ and a point).
\end{rem}

\begin{rem} Parts of the lemmas above can also be deduced from results in
\cite{McCord2}. For a related construction see also \cite{Bartsch}.
\end{rem}

We discuss now the following particular case. 
Consider a riemannian fibre bundle $\eta :\mathbf{R}^{n}\longrightarrow T\longrightarrow B$. 
As before let $\gamma$ be a flow on $B$ and let
$\gamma''$ be a lift of it to $T$ (of course, such lifts always exist). Assume that 
 there is a fixed quadratic form $q:\mathbf{R}^{n}\longrightarrow \mathbf{R}$ and 
that there exist local charts $U_{i}\subset B$ such that the restriction of $\eta$
to each $U_{i}$ is trivial and, moreover, with respect to these trivializations the vector 
field associated to $\gamma''$ splits as a direct sum of the vector field associated
to $\gamma$ and  $-\nabla q$; also assume that $q$ is compatible with the metric.

There is a particular locally trivial fibration that is associated to this context. 
Consider the associated sphere bundle of $\eta$:
$S(\eta)\stackrel{p}{\longrightarrow} B$. Denote by $e(q)$ the locally trivial
fibration $e(q)\longrightarrow B$ whose total space is the subset of 
$S(\eta)$ where $q$ is negative
or null and whose projection is the restriction of $p$. On $e(q)$
there is an obvious lift of $\gamma$ obtained by projecting first $\gamma''$ on
$S(\eta)$ and then by restriction.

\begin{prop} Let $S$ be an isolated invariant set of $\gamma$.
The Thom index $\overline{c}_{\gamma}^{e(q)}(S)$ is isomorphic to
$c_{\gamma''}(S)$ and this
identification is natural with respect to attractor-repellor pairs.
More precisely, if $(S;A,A^{\ast})$ is an attractor-repellor decomposition, then
$\overline{c}_{\gamma}^{e(q)}(S;A,A^{\ast})$ is isomorphic to 
$c_{\gamma''}(S;A,A^{\ast})$
\end{prop}

\begin{proof} Suppose $(N_{1},N_{0})$ is an index pair for $S$ with respect to
the flow $\gamma$. Let $N_{i}'$ be the total space of the restriction of
the unit disk bundle associated to $\eta$ restricted to $N_{i}$, $i=0,1$. 
 Let $e(N_{i})$ be the total space of the restriction of $e(q)$ to $N_{i}$. 
It is easy to see that the pair $(N_{1}', N_{0}'\bigcup e(N_{1}))$
is an index pair for $S$ with respect to the flow $\gamma''$.
Notice that if $D^{n}$ is the unit disk in $\mathbf{R}^{n}$ and $S^{n-1}$ its boundary, then
the region where $q$ is negative or null in $D^{n}$ has the structure of a cone
over the boundary. This implies that the region 
$L \subset N_{1}'$ where $q$ is negative
or null is homeomorphic to $N_{1}\bigcup (e(N_{1})\times [0,1])$ the
identification used making $x\times {0}$ correspond to the point $p(x)$ for 
each $x\in e(N_{1})$.
 On the other
hand we have a map of pairs 
$(L, (N_{0}'\bigcap L)\bigcup e(N_{1}))\subset (N_{1}', N_{0}'\bigcup e(N_{1})$.
This induces a map between the respective quotients:
$t:L/((N_{0}'\bigcap L)\bigcup e(N_{1}))\longrightarrow N_{1}'/(N_{0}'\bigcup 
e(N_{1}))$.
Because of the description of $L$ above, 
the domain of this map is a space in $\overline{c}_{\gamma}^{e(q)}$,
the target is a space in $c_{\gamma''}(S)$. This map is a homotopy equivalence
as both inclusions $L\subset N_{1}'$ and
$(N_{0}'\bigcap L)\bigcup e(N_{1})\subset N_{0}'\bigcup e(N_{1})$ are homotopy
equivalences (the first because the domain and the target
 are homotopy equivalent to $N_{1}$ and the second for a similar reason).

It is easy to verify, as in the lemma above, that the map $t$ induces an isomorphism
of connected simple systems and that this morphism is natural with respect to 
attractor-repellor pairs.
\end{proof}

\begin{rem} a. The proposition implies that the Thom index and the
corresponding cofibration sequence are completely determined by the
homotopical data extracted from the flow $\gamma''$.  

b. It is easy to see that the statement and proof above can be extended
in many ways. For example,
to the case when $q$ is a more complicated germ of an isolated but
still "reasonable" singularity (see \cite{Cornea1} for a definition of such 
reasonable critical points).
\end{rem}

\begin{cor} In the context of the proposition above,  if $\eta$ is trivial in a 
neighborhood of $S$, then we have
$c_{\gamma''}(S)\simeq\overline{c}_{\gamma}^{e(q)}(S)\simeq\Sigma^{k}c_{\gamma}(S)$
and if $(S;A,A^{\ast})$ is an attractor-repellor decomposition\\
 $c_{\gamma''}(S;A,A^{\ast})\simeq\overline{c}_{\gamma}^{e(q)}(S;A,A^{\ast})\simeq\Sigma^{k}c_{\gamma}(S;A,A^{\ast})$
where $k$ is the index of the quadratic form $q$.
\end{cor}

\begin{proof} Let $K=\{x\in D^{n} : q(x)\leq 0\}$ and let $H=K\bigcap \partial D^{n}$.
Notice that $H\he S^{k-1}$ and $K$ is contractible.
With the notations in the proof of the proposition, we have $L= N_{1}\times K$,
$N'_{0}\bigcap L=N_{0}\times K$, $e(N_{i})=N_{i}\times H$.
The statement is an immediate consequence of the construction of the Thom index
and of the fact that the cofibre of the projection $A\times S^{k-1}\longrightarrow A$
is, in a natural way, homotopy equivalent to $\Sigma^{k}A\vee S^{k}$. With this 
identification we have the commutative diagram below:
$$ \xy\xymatrix@-2pt{
N_{0}\times H\ar[d]\ar[r]&N_{1}\times H\ar[d] \\
N_{0}\times K\ar[r]\ar[d]&N_{1}\times K\ar[d]\\
\Sigma^{k} N_{0}\vee S^{k}\ar[r]&\Sigma^{k}N_{1}\vee S^{k}}
\endxy $$
the columns being (homotopy) cofibration sequences. By definition 
the Thom index of $S$ is
represented by the cofibre of the bottom horizontal map which is 
$\Sigma^{k} (N_{1}/N_{0})$.
\end{proof}

\begin{rem} The above corollary is well-known at the space level (in that case being
 a consequence of the product formula for the Conley index \cite{Conley}). 
For our purposes the key consequence of this corollary is the description 
of the connection map of $\gamma''$ in terms of that of $\gamma$. In particular,
when the fibration is trivial they are related by suspension.
\end{rem}

\subsection{Duality.}
Assume now that $B$ is a smooth, compact manifold of dimension $n$ that is 
embedded in a
sphere $S^{n+k}$. Let $\eta: \mathbf{R}^{k}\longrightarrow T\longrightarrow B$ 
be its normal bundle. Assume a riemannian metric
fixed on the total space of $\eta$. Let $S(\eta)$ be the induced unit sphere bundle.
We recall that we denote by $-\gamma$ the inverse flow of $\gamma$:
$-\gamma_{t}(x)=\gamma_{-t}(x)$. 
 
\begin{theo}Let $(S;A,A^{\ast})$ be an attractor-repeller pair for $\gamma$.
The cofibration sequences $\overline{c}_{\gamma}^{S(\eta)}(S;A,A^{\ast})$
and $c_{-\gamma}(S;A^{\ast},A)$ are Spanier-\\Whitehead duals by a duality
map that depends only on the embedding of $B$ in $S^{n+k}$.
\end{theo}

\begin{proof} 
We will construct a duality between two cofibration sequences
representing those of the statement. Afterwards we show 
 that the induced
duality of cofibrations of connected simple systems does not depend on the
choice of the representing sequences. The proof has four steps.

\subsubsection*{Reduction to a flow on $S^{n+k}$.}
Let $U$ be a tubular neighborhood of $B$. In the following we
identify it to $T$. Let $q: B\longrightarrow \mathbf{R}$ be a smooth function which
locally has the form $q(x,y)=-||x||^{2}$ where $x$ is the coordinate along the fibre of
$\eta$ and $y$ the coordinate along the base. Consider on $U$
the flow $\gamma''$ induced by the sum of the canonical horizontal lift of
$\gamma$ and of the flow induced by $-\nabla q$. 
Notice that $e(q)=S(\eta)$ and $e(-q)=\emptyset$. Applying Proposition 3.3 to $\gamma''$
 we obtain isomorphisms
$\overline{c}_{\gamma}^{S(\eta)}(S;A,A^{\ast})\he c_{\gamma''}(S;A ,A^{\ast})$,
$c_{-\gamma''}(S;A^{\ast},A)\he \overline{c}_{-\gamma}(S;A^{\ast},A)\he 
c_{-\gamma}(S;A^{\ast},A)$. This reduces the proof to showing that
$c_{\gamma''}(S;A,A^{\ast})$ and $c_{-\gamma''}(S;A^{\ast},A)$ are
Spanier-Whitehead duals.

\subsubsection*{Construction of some special index pairs.}
Extend $\gamma''$ to a flow on
$S^{n+k}$ which is constant outside a neighborhood of $U$. We continue to
 denote this flow by $\gamma''$. Let $S_{1}$ be the maximal attractor
of $\gamma''$ that does not contain $S$ , let $R_{1}$ be the complementary
repellor. Let $S_{2}$ be the attractor of $\gamma''$ obtained as the union
of $S_{1}$ with the set of the points of $A$ or situated on flow lines originating in $A$
and let $R_{2}$ be the complementary
repellor. Similarly, let $S_{3}$ be given by the union of $S_{2}$ with the set of the points
of $A^{\ast}$ or situated on flow lines originating at $A^{\ast}$, let $R_{3}$ 
be the complementary repellor. For each
 pair $(S_{i},R_{i})$ for $1\leq i\leq3$  we may construct 
by classical techniques (see, for example, \cite{Akin})
 smooth Lyapounov functions, 
$f_{i} :S^{n+k}\longrightarrow \mathbf{R}$  such that $f_{i}(S_{i})=1$ and 
$f_{i}(R_{i})=0$. By Sard's theorem we may find regular values $a$ and $b$ of
respectively $f_{1}$ and $f_{3}$ such that $f_{1}^{-1}(a)$ and $f_{3}^{-1}(b)$
intersect transversely. Notice that $S$ is the maximal invariant set inside
$N_{2}=\{ x : f_{1}(x)\geq a, f_{3}(x)\leq b\}$. Moreover let $U_{0}=
f_{1}^{-1}(a)\bigcap N_{2}$ and $V_{0}=f_{3}^{-1}(b)\bigcap N_{2}$. Then
$\partial N_{2}=U_{0}\bigcup V_{0}$ and the pairs $(N_{2},U_{0})$ and
$(N_{2},V_{0})$ are index pairs for $S$ for, respectively,
the flows $\gamma''$ and $-\gamma''$. We now choose a regular value $c$
of $f_{2}$ such that $f_{2}^{-1}(c)$ cuts transversely $f_{3}^{-1}(b)$ and
such that the inequality $f_{1}(x)\leq a$ implies $f_{2}(x) < c$ (this is possible
because $R_{2}\subset R_{1}$). Let $U_{1}=f_{2}^{-1}(c)\bigcap N_{2}$,
 $N_{1}=f_{2}^{-1}((-\infty,c])\bigcap N_{2}$, 
$L_{1}=N_{2}-(Int(N_{1})\bigcup U_{0})$, $G_{1}=\overline{N_{2}-N_{1}}$
and $V_{1}=\partial G_{1}\bigcap V_{0}$, $U_{2}=U_{1}\bigcup (V_{0}-V_{1})$. 
With these notations
the pairs $(N_{2},N_{1})$, $(L_{1},U_{2})$ and $(G_{1},U_{1})$ are index pairs of 
$A^{\ast}$ for the flow $\gamma''$; $(L_{1},V_{0})$ and $(G_{1}, V_{1})$
are index pairs for $A^{\ast}$ with respect to the flow $-\gamma''$;
$(N_{1},U_{0})$ is an index pair for $A$ with respect to $\gamma''$;
$(N_{2}, L_{1})$, $(N_{1}, U_{2})$  are index pairs for $A$ with respect to $-\gamma''$.

\subsubsection*{Duality.}
The basic argument that will be used is very classical in nature and is a variant of one that
appears in \cite{Atiyah}. For completeness we will formulate it as a lemma and
we will indicate the idea of the proof.
\begin{lem}
Assume that $M$ is an $n+k$ manifold with boundary
emdedded in $S^{n+k}$ such that $\partial M$ admits a decomposition $\partial M=
M_{0}\bigcup M_{1}$ with $M_{0}$ and $M_{1}$ being $n+k-1$-dimensional manifolds
such that $M_{1}\bigcap M_{0}=\partial M_{1}=\partial M_{0}$.
Then $M/M_{0}$ and $M/M_{1}$ are $n+k$-Spanier-Whitehead duals.
\end{lem}

\textbf{Proof of the lemma.} Embedd $S=S^{n+k}$ in $S^{n+k+1}$ as the boundary
of a disk $D=D^{n+k+1}$. Denote by $D'$ the complementary disk. 
Consider a copy of $M$, $M'$, in the exterior of $D$
and paste it, in $S^{n+k+1}$, to $D$ following $M_{0}$. 
(We use here and below the existence of collared neighborhoods
of $M$, $M_{0}$, $M_{1}$ and $M_{0}\bigcap M_{1}$ .) The result
of this pasting, $M^{\ast}$, is homotopy equivalent to $M/M_{0}$. 
As  $D'$ is the cone over the boundary of $D$ it is easy to see that
$M''=S^{n+k+1}-M^{\ast}$ deforms to the union  $M'''\bigcup_{M_{1}} D''$  with
 $D''$ a slightly smaller disk than $D'$ and $M'''$ a copy of $M$ that does not touch
$M'$. Hence, $M''\simeq M/M_{1}$.

\

Notice that the argument of the lemma can be also used to prove the duality of some maps.
Indeed assume that $K$ is an $n+k$-dimensional submanifold of $M$ such that $\partial K$
decomposes, similarly to $\partial M$, as $K_{0}\bigcup K_{1}$ and such that
  $K_{0}\subset M_{0}$, $K_{1}$  separates and is an NDR in $M$, $K_{0}$ and 
$K_{1}\bigcap \partial M$
are NDRs in $\partial M$. Then the maps
$K/K_{0}\longrightarrow M/M_{0}$ and $M/M_{1}\longrightarrow K/K_{1}$
are Spanier-Whitehead duals. The reason is that, as $K_{0}\subset M_{0}$
and $K/K_{1}=M/((M-K)\bigcup K_{1})$,  by the construction of the lemma
the two maps are identified to a pair of dual inclusions.\medskip

We return now to the proof of the theorem. We intend to show that the
cofibration sequences: $N_{1}/U_{0}\stackrel{i}{\longrightarrow}N_{2}/U_{0}
\stackrel{p}{\longrightarrow}N_{2}/N_{1}$ and
$L_{1}/V_{0}\stackrel{i'}{\longrightarrow}N_{2}/V_{0}\stackrel{p'}{\longrightarrow}
N_{2}/L_{1}$ are Spanier-Whitehead duals. Consider first the pair of maps $(i,p')$.
We have $N_{2}/L_{1}=N_{1}/U_{2}$ and the duality is immediate by the above.
Similarly, for the pair $(p,i')$ we have $L_{1}/V_{0}=G_{1}/V_{1}$
and $N_{2}/N_{1}=G_{1}/U_{1}$. 

\subsubsection*{Duality for connected simple systems.}
The proof is completed by showing that
 the duality induced by the two cofibration sequences described above
 does not depend on the choice of the
Lyapounov functions or on the choice of the constants $a$, $b$ and $c$.
It is easely seen that this statement reduces to proving that if
$N'_{2}$, $U'_{0}$ and $V'_{0}$ are constructed in the same way as,
respectively, $N_{2}$, $U_{0}$ and $V_{0}$ and we also have
$N'_{2}\subset Int (N_{2})$, then the map
$N'_{2}/U'_{0}\wedge N'_{2}/V'_{0}\stackrel{\mu'}{\longrightarrow} S^{n+k}$
is homotopic to the composition 
$N'_{2}/U'_{0}\wedge N'_{2}/V'_{0}\stackrel{c\wedge c^{\ast}}
{\longrightarrow}N_{2}/U_{0}\wedge N_{2}/V_{0}
\stackrel{\mu}{\longrightarrow}S^{n+k}$

Here $\mu$ and $\mu'$ are duality maps (here and above all duality maps
given by the lemma are constructed with respect to the same 
embedding of  $S^{n+k}$ in $S^{n+k+1}$); the maps
$c:N'_{2}/U'_{0}\longrightarrow N_{2}/U_{0}$ and
$c^{\ast}:N'_{2}/V'_{0}\longrightarrow N_{2}/V_{0}$  are the standard comparison
maps given, respectively, by  the connected simple systems $c_{\gamma''}(S)$ and
$c_{-\gamma''}(S)$. 

This can be seen as follows. 
Let $N''_{2}=\{x \in N_{2} : \exists t\geq 0, \gamma''_{-t}(x)\in N'_{2}\}$,
 $U''_{0}=N''_{2}\bigcap U_{0}$ and  
$V''_{0}=\overline{\partial N''_{2}-U''_{0}}$. Notice that $U'_{0}$
and $U''_{0}$ are homeomorhpic as well as the pairs $V'_{0}$,
$V''_{0}$ and $N'_{2}$, $N''_{2}$. The comparison map $c$
factors using these homeomorphisms 
$N'_{2}/U'_{0}\longrightarrow N''_{2}/U''_{0}\stackrel{c''}{\longrightarrow}
 N_{2}/U_{0}$ with $c''$ being induced by the inclusions.
We apply the general duality argument in the lemma to obtain that
$c''$ is dual to the map $N_{2}/V_{0}\longrightarrow N''_{2}/V''_{0}$.
This map is clearly homotopic to the comparison map $c''':N_{2}/V_{0}
\longrightarrow N''_{2}/V''_{0}$.

Hence, if we denote by $c^{\ast}_{1}$ the comparison map 
which is the homotopy inverse of
$c^{\ast}$ we obtain that $c^{\ast}_{1}$ and $c$ are duals and this implies
the homotopy we look for. 
 \end{proof} 

\begin{cor}
The two connection maps $\overline{\delta}:\overline{c}^{S(\eta)}_{\gamma}(A^{\ast})
\longrightarrow 
\Sigma \overline{c}^{S(\eta)}_{\gamma}(A)$ and $\delta^{\ast}: c_{-\gamma}(A)
\longrightarrow \Sigma c_{-\gamma}(A^{\ast})$ are Spanier-Whitehead duals. In particular,
if $\eta$ is trivial, then $\delta^{\ast}$ is dual to the connection map
$\delta: c_{\gamma}(A^{\ast})\longrightarrow \Sigma c_{\gamma}(A)$.
\end{cor}

\begin{rem}
This clearly extends the relation, discovered by Franks \cite{Franks}, 
between the two connection maps in the
case when $A$ and $A^{\ast}$ are two succesive critical points of a Morse-Smale flow.
The corollary also implies, via the Thom isomorphism, the (co)-homological
results of McCord \cite{McCord}. 
It also provides an extension of the
duality result for succesive reasonable critical points of \cite{Cornea2}. However,
in the Morse-Smale case the understanding of the geometry of the problem
is better because of the explicit description of the stable and unstable manifolds
of the critical points. Also, more precise information can be deduced in 
the  case of succesive reasonable critical points
with $\eta$  non trivial.  
In that case $\eta$ is trivial in  small enough neighborhoods of $A$ and of $A^{\ast}$,
hence $\Sigma^{k}c_{\gamma}(A)\he \overline{c}_{-\gamma}^{S(\eta)}(A)$
and similarly for $A^{\ast}$. One can then prove that 
$\delta$ and $\delta^{\ast}$ are Spanier-Whitehead duals \textit{modulo some
twisting} depending in a precise way on $\eta$.

\end{rem}

\bibliographystyle{amsplain}

\end{document}